\documentclass[onecolumn]{autart}
\usepackage{graphicx}
\input{amssym}
\def\be{\begin{eqnarray}}
\def\ee{\end{eqnarray}}

\def\di{\displaystyle }

\begin{document}
\begin{frontmatter}
\title{Lie group analysis of Poisson's equation \\ and optimal system of
subalgebras for \\ Lie algebra of $3-$dimensional rigid motions} 
\thanks[footnoteinfo]{Corresponding author: Tel. +98-21-73913426.
Fax +98-21-77240472.}
\author[Tehran]{M. Nadjafikhah \thanksref{footnoteinfo}}\ead{m\_nadjafikhah@iust.ac.ir}
\address[Tehran]{School of Mathematics, Iran University of Science and Technology, Narmak-16, Tehran, I.R.IRAN.}
\renewcommand{\sectionmark}[1]{}
\begin{abstract}
Using the basic Lie symmetry method, we find the most general Lie
point symmetries group of the $\nabla u=f(u)$ Poisson's equation,
which has a subalgebra isomorphic to the $3-$dimensional special
Euclidean group ${\rm SE}(3)$ or group of rigid motions of ${\Bbb
R}^3$. Looking the adjoint representation of ${\rm SE}(3)$ on its
Lie algebra $\goth{se}(3)$, we will find the complete optimal
system of its subalgebras. This latter provides some properties of
solutions for the Poisson's equation.
\end{abstract}
\end{frontmatter}
%
%
\section{Introduction}
The {\it Poisson's equation} is
\be \nabla u=f(u), \label{eq}\ee
where $\nabla={\partial^2}/{\partial x^2}+{\partial^2}/{\partial
y^2}+{\partial^2}/{\partial z^2}$, $f:{\Bbb R}\to{\Bbb R}$ is a
given smooth function, $\Omega$ is an open subset of ${\Bbb R}^3$
and  $u:\Omega\to{\Bbb R}$ is an unknown smooth function
\cite{Jost}. This is a stationary heat equation with a nonlinear
source (5.4.1 of \cite{PolZai}). In this paper, using the basic
Lie symmetry method, we find the most general Lie point symmetries
group of the Poisson's equation, which has a subalgebra isomorphic
to $3-$dimensional special Euclidean group ${\rm SE}(3)$. Looking
the adjoint representation of ${\rm SE}(3)$ on its Lie algebra
$\goth{se}(3)$, we will find the complete classification of its
subalgebras up to conjugacy. This latter provide a process for
building new solutions for the Poisson's equation.

Inspired by Galois' theory of, Sophus Lie developed an analogous
theory of symmetry for differential equations. Lie's theory led to
an algorithmic way to find special explicit solutions to
differential equations with symmetry. These special solutions are
called {\it group invariant solutions} and they constitute
practically every known explicit solution to the systems of
non-linear partial differential equations arising in mathematical
physics, differential geometry and other areas.

These group-invariant solutions are found by solving a reduced
system of differential equations involving fewer independent
variables than the original system. For example, the solutions to
a partial differential equation in two independent variables which
are invariant under a given one-parameter symmetry group are all
found by solving a system of ordinary differential equations.
Today the search for group invariant solutions is still a common
approach to explicitly solving non-linear partial differential
equations. A good modern introduction to the use of similarity
methods in engineering applications is the book \cite{SechNa}. An
excellent introductory reference to Lie's theory of group
invariant solutions is \cite{Olv1}.
\section{Lie point symmetries}
In this section , we find the Lie point symmetries of Poisson's
equation, using the basic prolongation method and the
infinitesimal criterion of invariance as \cite{Olv1}. Consider the
one parameter Lie group of infinitesimal transformations on
$(x_1=x,x_2=y,x_3=z,u)$ given by
\be \tilde{x}_i = x^i+s\,\xi_1(x_1,x_2,x_3,u)+O(s^2),\quad
\tilde{u} = u+s\,\varphi(x_1,x_2,x_3,u)+O(s^2),\quad i=1,2,3, \ee
where $s$ is the group parameter and $\xi_i$, $i=1,2,3$ and
$\varphi$ are the infinitesimals of the transformations for the
independent and dependent variables, respectively. The associated
vector field is ${\bf
v}=\sum_{i=1}^3\xi_i\partial_{x_i}+\varphi\,\partial_u$, and its
second prolongation is
\be
{\rm Pro}^{(2)}{\bf v}={\bf
v}+\varphi^x\,\partial_{u_x}+\varphi^y\,\partial_{u_y}+\varphi^z\,\partial_{u_z}+\varphi^{xx}\,\partial_{u_{xx}}+\varphi^{xy}\,\partial_{u_{xy}}
+\varphi^{xz}\,\partial_{u_{xz}}+\varphi^{yy}\,\partial_{u_{yy}}+\varphi^{yz}\,\partial_{u_{yz}}+\varphi^{zz}\,\partial_{u_{zz}},
\ee
with
\be
\varphi^x&=&-u_x(\partial_x+u_x\partial_u)\,\xi_1-u_y(\partial_x+u_x\partial_u)\,\xi_2
-u_z(\partial_x+u_x\partial_u)\,\xi_3+(\partial_x+u_x\partial_u)\,\varphi,\nonumber\\
&&\vdots\\
\varphi^{zz}&=&-(2u_{xz}\partial_z+(u_xu_{zz}+2u_zu_{xz})\partial_u-u_x\partial_{zz}-2u_xu_z\partial_{zu}-u_xu_z^2\partial_{uu})\,\xi_1\nonumber\\
&&-(2u_{yz}\partial_z+u_yu_{zz}+2u_zu_{yz})\partial_u+u_y\partial_{zz}+2u_yu_z\partial_{zu}+u_yu_z^2\partial_{uu})\,\xi_2\nonumber\\
&&-(2u_{zz}\partial_z+3u_{zz}u_z\partial_u+u_z\partial_{zz}+2u_z^2\partial_{zu}+u_z^3\partial_{uu})\,\xi_3
+(u_{zz}\partial_u+\partial_{zz}+2u_z\partial_{zu}+u_z^2\partial_{uu})\,\varphi.\nonumber
\ee
The vector field $\bf v$ generates a one parameter symmetry group
of Poisson's equation if and only if
\be {\rm Pro}^{(2)}{\bf v}(\nabla
u-f(u))=0\hspace{1cm}\mbox{whenever}\hspace{1cm} \nabla u=f(u).\ee
This condition is equivalent to the set of defining equations:
\be \begin{array}{l}
\partial_u\xi_1=0,\\
\partial_u\xi_2=0,\\
\partial_u\xi_3=0,\\
\partial_{uu}\varphi=0,
\end{array}\qquad
\begin{array}{l}
\partial_x\xi_1=\partial_z\xi_3,\\
\partial_x\xi_2+\partial_y\xi_1=0,\\
\partial_x\xi_3+\partial_z\xi_1=0,\\
\partial_y\xi_3+\partial_z\xi_2=0,\\
\partial_z\xi_3=\partial_y\xi_2,
\end{array}\qquad
\begin{array}{l}
\partial_{zz}\xi_2+\partial_{yy}\xi_2+\partial_{xx}\xi_2=2\partial_{yu}\varphi,\\
\partial_{zz}\xi_3+\partial_{yy}\xi_3+\partial_{xx}\xi_3=2\partial_{zu}\varphi,\\
\partial_{xx}\xi_1+\partial_{zz}\xi_1+\partial_{yy}\xi_1=2\partial_{xu}\varphi,\\
\end{array}\qquad \nabla\varphi=2F(u).\partial_z\xi_3+F'(u).\varphi.\ee
By solving this system of PDEs, we find that
\paragraph{Theorem}
{\em The Lie group $\goth{g}_f$ of point symmetries of Poisson's
equation (\ref{eq}) is an infinite dimensional Lie algebra
generated by the vector fields ${\bf
X}=\xi_1\partial_x+\xi_2\partial_y+\xi_3\partial_z+\varphi\partial_u$,
where
\be
\xi_1&=&a_7(x^2-y^2-z^2)+2(a_5y+a_1z)x+a_6x+a_8y-a_4z+a_9,\nonumber \\
\xi_2&=&a_5(y^2-z^2-x^2)+2(a_7x+a_1z)y-a_8x+a_6y+a_2z+a_{10},\qquad \varphi=F_1(x,y,z).u+F_2(x,y,z), \\
\xi_3&=&a_1(z^2-x^2-y^2)+2(a_7x+a_5y)z+a_4x-a_2y+a_6z+a_3,\nonumber
\ee
with
\be F_1=a_{11}-(a_7x+a_5y+a_1z),\quad \nabla
F_2=\big(5F_1-2a_6-4a_{11}\big).f(u)-\big(uF_1+F_2\big).f'(u), \ee
and $a_i$, $i=1,\cdots,11$ are arbitrary constants.}
\paragraph{Conclusion}
{\em The Lie group $\goth{g}_f$ of point symmetries of Poisson's
equation (\ref{eq}) has a $6-$dimensional Lie subalgebra generated
by the following vector fields:
\begin{eqnarray} \label{6vf}
X_1=\partial_x,\quad X_2=\partial_y,\quad X_3=\partial_z, \quad
X_4=y\partial_z-z\partial_y,\quad
X_5=z\partial_x-x\partial_z,\quad X_6=x\partial_y-y\partial_x.
\end{eqnarray}
The commutator table of $\goth{g}_f$ is given in Table \ref{ta-1},
where the entry in the $i^{\rm th}$ row and $j^{\rm th}$ column is
defined as $[X_i,X_j]=X_i.X_j-X_j.X_i$, $i,j=1,\cdots,8$.\hfill\
$\Box$
\begin{table}
\caption{Commutation table of $\goth{g}_f$.}\label{ta-1}
\begin{tabular}{ccccccc} \hline &$X_1$&$X_2$&$X_3$&$X_4$&$X_5$&$X_6$\\ \hline
$X_1$&$0   $&$0   $&$0   $&$0   $&$-X_3$&$X_2 $\\
$X_2$&$0   $&$0   $&$0   $&$X_3 $&$0   $&$-X_1$\\
$X_3$&$0   $&$0   $&$0   $&$-X_2$&$X_1 $&$0   $\\
$X_4$&$0   $&$-X_3$&$X_2 $&$0   $&$-X_6$&$X_5 $\\
$X_5$&$X_3 $&$0   $&$-X_1$&$X_6 $&$0   $&$-X_4$\\
$X_6$&$-X_2$&$X_1 $&$0   $&$-X_5$&$X_4 $&$0   $ \\ \hline
\end{tabular}
\end{table}}

An isometry (or rigid motion) of Euclidean space ${\Bbb R}^3$ is a
mapping ${\Bbb R}^3\to{\Bbb R}^3$ which preserves the Euclidean
distance between points. For every isometry $F$, there exist a
unique transformation $T\in{\rm Trans}(3)\cong{\Bbb R}^3$ and a
unique orthogonal transformation $C\in {\rm O}(3)$ such that
$F=C\circ T$. ${\rm Trans}(3)$ is a normal subgroup of ${\rm
E}(3)$, and ${\rm E}(3)\cong{\Bbb R}^3\propto{\rm O}(3)$. ${\rm
E}(3)$ has a natural $6-$dimensional Lie group structure. The
connected component containing $I_3$ is {\it special Euclidean
group} ${\rm SE}(3)$, consists of all orientation preserving
isometries, and is a connected $6-$dimensional Lie group.
Furthermore, the Lie algebra $\goth{se}(3)$ of $3-$dimensional
special Euclidean group is spanned by the following 6 vector
fields (\ref{6vf}).
\paragraph{Remark}
The ${\goth g}_f$ Lie group of point symmetries of Poisson's
equation $\nabla u=f(u)$ may be changed, for special choices of
$f$; for example, if $f=0$, that is for $\nabla u=0$ Laplace
equation, ${\goth g}_f$ is spanned by $X_i$, $i=1,\cdots,6$,
$\partial_u$, $\partial_x+\partial_y+\partial_z$,
$xz\partial_x+yz\partial_y+(z^2-x^2-y^2)\partial_z+zu\partial_u$,
$(x^2-y^2-z^2)\partial_x+xy\partial_y+xz\partial_z+xu\partial_u$,
$xy\partial_x+(x^2+y^2+z^2)\partial_y+yz\partial_z+yu\partial_u$,
and so on.
\section{Optimal system of subalgebras}
Let a system of differential equation $\Delta$ admitting the
symmetry Lie group $G$,be given. Now $G$ operates on the set of
solutions $S$ of $\Delta$. Let $s\cdot G$ be the orbit of $s$, and
$H$ be a subgroup of $G$. Invariant $H-$solutions $s\in S$ are
characterized by equality $s\cdot S=\{s\}$. If $h\in G$ is a
transformation and $s\in S$,then $h\cdot(s\cdot H)=(h\cdot s)\cdot
(hHh^{-1})$. Consequently,every invariant $H-$solution $s$
transforms into an invariant $hHh^{-1}-$solution (Proposition 3.6
of \cite{Olv1}).

Therefore, different invariant solutions are found from similar
subgroups of $G$. Thus,classification of invariant $H-$solutions
is reduced to the problem of classification of subgroups of $G$,up
to similarity. An optimal system of $s-$dimensional subgroups of
$G$ is a list of conjugacy inequivalent $s-$dimensional subgroups
of $G$ with the property that any other subgroup is conjugate to
precisely one subgroup in the list. Similarly, a list of
$s-$dimensional subalgebras forms an optimal system if every
$s-$dimensional subalgebra of $\goth g$ is equivalent to a unique
member of the list under some element of the adjoint
representation: $\tilde{\goth h}={\rm Ad}(g)\cdot{\goth h}$.

Let $H$ and $\tilde{H}$ be connected, $s-$dimensional Lie
subgroups of the Lie group $G$ with corresponding Lie subalgebras
${\goth h}$ and $\tilde{\goth h}$ of the Lie algebra ${\goth g}$
of $G$. Then $\tilde{H}=gHg^{-1}$ are conjugate subgroups if and
only $\tilde{\goth h}={\rm Ad}(g)\cdot{\goth h}$ are conjugate
subalgebras (Proposition 3.7 of \cite{Olv1}). Thus,the problem of
finding an optimal system of subgroups is equivalent to that of
finding an optimal system of subalgebras, and so we concentrate on
it.

For one-dimensional subalgebras,the classification problem is
essentially the same as the problem of classifying the orbits of
the adjoint representation,since each one-dimensional subalgebra
is determined by a nonzero vector in $\goth{se}(3)$ and so to
"simplify" it as much as possible.

The adjoint action is given by the Lie series
\be \mathrm{Ad}(\exp(s.X_i).X_j) =
X_j-s.[X_i,X_j]+\frac{s^2}{2}.[X_i, [X_i,X_j]]-\cdots,\ee
where $[X_i,X_j]$ is the commutator for the Lie algebra, $s$ is a
parameter, and $i,j=1,\cdots,6$ (\cite{Olv1},page 199). We can
write the adjoint action for the Lie algebra $\goth{se}(3)$, and
show that
\paragraph{Theorem}
{\em A one-dimensional optimal system of $\goth{se}(3)$ is given
by
\be
\begin{array}{l}
A_1^1 : X_6,\\[1mm]A_1^2 : X_1+b\,X_4,\end{array}\;\;\;\;\;
\begin{array}{l}
A_1^3 : X_2+b\,X_5,\\[1mm]A_1^4 :
X_3+b\,X_6,\end{array}\;\;\;\;\;
\begin{array}{l}
A_1^5 : X_1+a\,X_2+b\,X_6,\\[1mm]A_1^6 : X_1+a\,X_3+b\,X_4,\end{array}\;\;\;\;\;A_1^7
:X_2+a\,X_3+bX_5, \ee
where $a,b\in{\Bbb R}$ and $a\neq0$.}

\medskip \noindent {\it Proof:} $F^s_i:\goth{se}(3)\to \goth{se}(3)$ defined by
$X\mapsto\mathrm{Ad}(\exp(s_iX_i).X)$ is a linear map,for
$i=1,\cdots,6$. The matrix $M^s_i$ of $F^s_i$, $i=1,\cdots,6$,with
respect to basis $\{X_1,\cdots,X_6\}$ is
\be \label{eq-3}
&\di M_1^s=\left[ \begin {array}{cccccc}
1&0&0&0&0&0\\0&1&0&0&0&0\\0&0&1&0&0&0\\0&0&0&1&0&0\\0&0&s&0&1&0\\0&-s&0&0&0&1
\end{array} \right], \qquad
M_2^s=\left[ \begin {array}{cccccc}
1&0&0&0&0&0\\0&1&0&0&0&0\\0&0&1&0&0&0\\0&0&-s&1&0&0\\0&0&0&0&1&0\\s&0&0&0&0&1
\end{array} \right],\qquad
M_3^s=\left[ \begin {array}{cccccc}
1&0&0&0&0&0\\0&1&0&0&0&0\\0&0&1&0&0&0\\0&s&0&1&0&0\\-s&0&0&0&1&0\\0&0&0&0&0&1
\end{array} \right],\nonumber&\\[-2mm] \\[-2mm]
&\di M_4^s=\left[ \begin {array}{cccccc}
1&0&0&0&0&0\\0&C&S&0&0&0\\0&S&C&0&0&0\\0&0&0&1&0&0\\0&0&0&0&C&S\\0&0&0&0&-S&C
\end {array} \right],\qquad
M_5^s=\left[ \begin {array}{cccccc}
C&0&-S&0&0&0\\0&1&0&0&0&0\\S&0&C&0&0&0\\0&0&0&C&0&-S\\0&0&0&0&1&0\\0&0&0&S&0&C
\end{array} \right],\qquad
M_6^s=\left[ \begin {array}{cccccc}
C&S&0&0&0&0\\-S&C&0&0&0&0\\0&0&1&0&0&0\\0&0&0&C&S&0\\0&0&0&-S&C&0\\0&0&0&0&0&1
\end{array} \right].\nonumber&
\ee
respectively,where $S=\sin s$ and $C=\cos s$. Let
$X=\sum_{i=1}^6a_iX_i$,then
\be
F^{s_6}_6\circ\cdots\circ F^{s_1}_1\;:\;X\;\mapsto\; (\cos
s_5.\cos s_6.a_1+\cos s_5.\sin s_6.a_2-\sin s_5.a_3).X_1+\cdots
\ee
Now,we can simply $X$ as follows:

If $a_1$, $a_2$ and $a_3=0$,then we can make the coefficients of
$X_1,\cdots,X_5$ vanish,by $s_4=-\arctan(a_5/a_6)$ and
$s_5=\arctan(a_4/a_6)$. Scaling $X$ if necessary,we can assume
that $a_6=1$. And $X$ is reduced to the Case of $A_1^1$.

If $a_2$ and $a_3=0$ but $a_1\neq0$,then we can make the
coefficients of $X_2$, $X_3$, $X_5$ and $X_6$ vanish,by
$s_2=-a_6/a_1$ and $s_3=a_5/a_1$. Scaling $X$ if necessary,we can
assume that $a_1=1$. And $X$ is reduced to the Case of $A_1^2$.

If $a_1$ and $a_3=0$ but $a_2\neq0$,then we can make the
coefficients of $X_1$, $X_2$, $X_3$, and $X_6$ vanish,by
$s_1=a_6/a_2$ and $s_3=-a_4/a_2$. Scaling $X$ if necessary,we can
assume that $a_3=1$. And $X$ is reduced to the Case of $A_1^3$.

If $a_1$ and $a_2=0$ but $a_3\neq0$,then we can make the
coefficients of $X_1$, $X_2$, $X_4$ and $X_5$ vanish,by
$s_1=-a_5/a_3$ and $s_2=a_4/a_3$. Scaling $X$ if necessary,we can
assume that $a_3=1$. And $X$ is reduced to the Case of $A_1^4$.

If $a_1$ and $a_2\neq0$,then we can make the coefficients of
$X_3$, $X_4$ and $X_5$ vanish,by $s_1=a_6/a_2$, $s_3=a_5/a_1$ and
$s_4=-\arctan(a_3/a_2)$. Scaling $X$ if necessary,we can assume
that $a_1=1$. And $X$ is reduced to the Case of $A_1^5$.

If $a_2=0$ but $a_1$ and $a_3\neq0$,then we can make the
coefficients of $X_2$, $X_5$ and $X_6$ vanish,by $s_1=-a_5/a_3$,
$s_2=-a_6/a_1$ and $s_4=-\arctan(a_3/a_2)$. Scaling $X$ if
necessary,we can assume that $a_1=1$. And $X$ is reduced to the
Case of $A_1^6$.

If $a_1=0$ but $a_2$ and $a_3\neq0$,then we can make the
coefficients of $X_3$, $X_4$ and $X_5$ vanish,by $s_1=a_6/a_2$,
and $s_2=a_4/a_3$. Scaling $X$ if necessary,we can assume that
$a_2=1$. And $X$ is reduced to the Case of $A_1^7$. \hfill\ $\Box$

\paragraph{Theorem}
{\em A two-dimensional optimal system of $\goth{se}(3)$ is given
by
\be
\begin{array}{l}
A_2^1 : X_2,X_5,\\[1mm]
A_2^2 : X_3,X_6,\end{array}\;\;\;\;\
\begin{array}{l}
A_2^3 : X_1,X_2+a\,X_4,\\[1mm]
A_2^4 : X_1,X_3+a\,X_4,\end{array}\;\;\;\;\
\begin{array}{l}
A_2^5 : X_2,X_3+a\,X_5,\\[1mm]
A_2^6 : X_3,X_1+a\,X_4,\end{array}
\ee
where $a\in{\Bbb R}$. All of these sub-algebras are Abelian.}

\medskip \noindent {\it Proof:} In this proof,we shall assume that
${\goth g}={\rm Span}_{\Bbb R}\{X,Y\}$ is a 2-dimensional Lie
sub-algebra of $\goth{se}(3)$.

Let $X=X_6$ be as the Case of $A_1^1$, $Y=\sum_{i=1}^6b_i\,X_i$,
and $[X,Y]=\alpha\,X+\beta\,Y$. Then,we have $Y=b_3\,X_3+b_6\,X_6$
and $\alpha=\beta=0$. By a suitable change of base of $\goth g$,we
can assume that $Y=X_3$. Now $\goth g$ is reduced to the Case of
$A_2^2$. We can not be used to further simplify this
sub-algebra,by $F^{s_i}_i$, $i=1,\cdots,6$ defined as
(\ref{eq-3}).

Let $X=X_1+a_6\,X_6$ be as the Case of $A_1^2$,
$Y=\sum_{i=1}^6b_i\,X_i$, and $[X,Y]=\alpha\,X+\beta\,Y$. Then,we
have $Y=b_3\,X_3+b_1\,X$ and $\alpha=\beta=0$. By a suitable
change of base of $\goth g$,we can assume that $Y=X_3$. Now $\goth
g$ is reduced to the Case of $A_2^6$. We can not be used to
further simplify this sub-algebra,by $F^{s_i}_i$, $i=1,\cdots,6$.

Let $X=X_2+a_5\,X_5$ be as the Case of $A_1^3$,
$Y=\sum_{i=1}^6b_i\,X_i$, and $[X,Y]=\alpha\,X+\beta\,Y$. Then,we
have $Y=b_2\,X_2+b_5\,X_5$ and $\alpha=\beta=0$. By a suitable
change of base of $\goth g$,we can assume that $X=X_2$ and
$Y=X_5$. Now, $\goth g$ is reduced to the Case of $A_2^1$. We can
not be used to further simplify this sub-algebra,by $F^{s_i}_i$,
$i=1,\cdots,6$.

Let $X=X_3+a_6\,X_6$ be as the Case of $A_1^4$,
$Y=\sum_{i=1}^6b_i\,X_i$, and $[X,Y]=\alpha\,X+\beta\,Y$. Then,we
have $Y=b_3\,X_3+b_6\,X_6$ and $\alpha=\beta=0$. By a suitable
change of base of $\goth g$,we can assume that $X=X_6$ and
$Y=X_3$. Now, $\goth g$ is reduced to the Case of $A_2^2$.

Let $X=X_1+a_2\,X_2+a_4\,X_4$ be as the Case of $A_1^5$,
$Y=\sum_{i=1}^6b_i\,X_i$, and $[X,Y]=\alpha\,X+\beta\,Y$. Then,we
have $Y=b_1\,X_1+(b_2/a_2)\,X$ and $\alpha=\beta=0$. By a suitable
change of base of $\goth g$,we can assume that $X=X_1$ and
$Y=X_2+a\,X_4$. Now, $\goth g$ is reduced to the Case of $A_2^3$.
We can not be used to further simplify this sub-algebra,by
$F^{s_i}_i$, $i=1,\cdots,6$.

Let $X=X_2+a_3\,X_3+a_4\,X_4$ be as the Case of $A_1^6$,
$Y=\sum_{i=1}^6b_i\,X_i$, and $[X,Y]=\alpha\,X+\beta\,Y$. Then,we
have $Y=b_1\,X_1+(b_3/a_3)\,X$ and $\alpha=\beta=0$. By a suitable
change of base of $\goth g$,we can assume that $X=X_1$ and
$Y=X_3+a\,X_4$. Now, $\goth g$ is reduced to the Case of $A_2^4$.
We can not be used to further simplify this sub-algebra,by
$F^{s_i}_i$, $i=1,\cdots,6$.

Let $X=X_2+a_3\,X_3+a_5\,X_5$ be as the Case of $A_1^7$,
$Y=\sum_{i=1}^6b_i\,X_i$, and $[X,Y]=\alpha\,X+\beta\,Y$. Then,we
have $Y=b_2\,X_2+(b_3/a_3)\,X$ and $\alpha=\beta=0$. By a suitable
change of base of $\goth g$,we can assume that $X=X_2$ and
$Y=X_3+a\,X_5$. Now, $\goth g$ is reduced to the Case of $A_2^5$.
We can not be used to further simplify this sub-algebra,by
$F^{s_i}_i$, $i=1,\cdots,6$. And,we have proved the Theorem.
\hfill\ $\Box$
\paragraph{Theorem}
{\em A three-dimensional optimal system of $\goth{se}(3)$ is given
by $A_3 : X_1+a\,X_4,X_2,X_3$, where $a\in{\Bbb R}$. The
commutator table of $A_3$ is Table \ref{ta-2}.
\begin{table}
\caption{Commutation table of $A_3$.}\label{ta-2}
\begin{tabular}{cccc} \hline &$X$&$Y$&$Z$\\\hline
&&& \\[-5mm] $X$&$0$   &$-a$  &$0$\\ Y&$a$   &$0$   &$0$\\ $Z$&$0$   &$0$   &$0$\\\hline
\end{tabular}
\end{table}}

\medskip \noindent {\it Proof:} In this proof,we assume that
${\goth g}={\rm Span}_{\Bbb R}\{X,Y,Z\}$ is a 3-dimensional Lie
sub-algebra of $\goth{se}(3)$.

Let $X=X_2$ and $Y=X_5$ are as the Case of $A_2^1$,
$Z=\sum_{i=1}^6b_i\,X_i$,
$[X,Z]=\alpha_1\,X+\beta_2\,Y+\beta_3\,Z$ and
$[X,Z]=\beta_1\,X+\beta_2\,Y+\beta_3\,Z$. Then,we have
$Z=b_2\,X+b_5\,Y$. By a suitable change of base of $\goth g$,we
can assume that $Z=0$, and $\goth g$ is not a $3-$dimensional
sub-algebra. Thus,in this case we have not any $3-$dimensional
sub-algebra.

The Cases $A_2^i$, $i=2,3,4,5$ are similar. Thus,in these cases we
have not any $3-$dimensional sub-algebra.

Let $X=X_3$ and $Y=X_1+a\,X_4$ are as the Case of $A_2^6$,
$Z=\sum_{i=1}^6b_i\,X_i$,
$[X,Z]=\alpha_1\,X+\beta_2\,Y+\beta_3\,Z$ and
$[X,Z]=\beta_1\,X+\beta_2\,Y+\beta_3\,Z$. Then,we have
$Z=(b_1/\alpha_3).\big(-\beta_3\,X+\alpha_3\,Y-a\,X_2\big)$. By a
suitable change of base of $\goth g$,we can assume that
$X=X_1+a\,X_4$, $Y=X_2$ and  $Z=X_3$. Now, $\goth g$ is reduced to
the Case of $A_3$. We can not be used to further simplify this
sub-algebra,by $F^{s_i}_i$, $i=1,\cdots,6$ defined as
(\ref{eq-3}). \hfill\ $\Box$
\paragraph{Theorem}
{\em A four-dimensional optimal system of $\goth{se}(3)$ is given
by $A_4 : X_1,X_2,X_3,X_4$. The commutator table of $A_4$ is Table
\ref{ta-3}.
\begin{table}
\caption{Commutation table of $A_4$.}\label{ta-3}
\begin{tabular}{ccccc} \hline&$X_1$&$X_2$&$X_3$&$X_4$\\\hline
$X_1$&$0$   &$0$   &$0$   &$0$   \\ $X_2$&$0$   &$0$   &$0$   &$X_3$ \\
$X_3$&$0$ &$0$ &$0$ &$-X_2$\\ $X_4$&$0$   &$-X_3$&$X_2$
&$0$\\\hline
\end{tabular}
\end{table}}
\medskip \noindent {\it Proof:} Assume that
${\goth g}={\rm Span}_{\Bbb R}\{X_1+a\,X_4,X_2,X_3,X\}$ be a Lie
sub-algebra of $\goth{se}(3)$,where $X=\sum_{i=1}^6b_i\,X_i$.
Then,we have $X=\sum_{i=1}^4b_i\,X_i$. By a suitable change of
base of $\goth g$,we can assume that $X=X_4$. Now, $\goth g$ is
reduced to the Case of $A_4$. We can not be used to further
simplify this sub-algebra,by $F^{s_i}_i$, $i=1,\cdots,6$ defined
as (\ref{eq-3}). \hfill\ $\Box$
\paragraph{Theorem}
{\em $\goth{se}(3)$ has not any five-dimensional Lie sub-algebra.}

\medskip \noindent {\it Proof:} Assume that
${\goth g}={\rm Span}_{\Bbb R}\{X_1,X_2,X_3,X_4,X\}$ be a Lie
sub-algebra of $\goth{se}(3)$,where $X=\sum_{i=1}^6b_i\,X_i$.
Then,we have $X=\sum_{i=1}^4b_i\,X_i$, and ${\goth g}={\rm
Span}_{\Bbb R}\{X_1,X_2,X_3,X_4\}$ is not a $5-$dimensional
sub-algebra. \hfill\ $\Box$
\section{Invariant solutions}
To obtain the group transformation which is generated by the
infinitesimal generators
$X_i=\sum_{j=1}^3\xi^i_j+\phi^i\partial_u$ for $i=1\cdots,6$, we
need to solve the 6 systems of first order ordinary differential
equations,
\be
\begin{array}{lcl}
\tilde{x}_j'(s) =
\xi^i_j(\tilde{x}_1(s),\tilde{x}_2(s),\tilde{x}_3(s),\tilde{u}(s)),
&& \tilde{x}_j(0)=x_j, \\[2mm]
\tilde{u}'(s) =
\phi^i(\tilde{x}_1(s),\tilde{x}_2(s),\tilde{x}_3(s),\tilde{u}(s)),
&& \tilde{u}(0)=u.
\end{array}\;\;\; \begin{array}{l} i=1,\cdots,6 \\ j=1,2,3
\end{array}
\ee
Exponentiating the infinitesimal symmetries of the Poisson's
equation (\ref{eq}), we get the one parameter groups $g_k(s)$
generated by $X_k$ for $k=1,\cdots,6$; and consequently, we have
\paragraph{Theorem}
{\em If $u=h(x,y,z)$ is a solution of the Poisson's equation
(\ref{eq}), so are the functions
\be
\begin{array}{l}
g_1(s)\cdot h(x,y,z)=h(x+s,y,z),\\
g_2(s)\cdot h(x,y,z)=h(x,y+s,z),\\
g_3(s)\cdot h(x,y,z)=h(x,y,z+s),\end{array}\qquad
\begin{array}{l}
g_4(s)\cdot h(x,y,z)=h(x,y\cos s-z\sin s,z\cos s+y\sin s),\\
g_5(s)\cdot h(x,y,z)=h(x\cos s+z\sin s,y,z\cos s-x\sin s),\\
g_6(s)\cdot h(x,y,z)=h(x\cos s-y\sin s,x\sin s+y\cos
s,z),\end{array}
\ee}

\begin{thebibliography}{99}%
%
\bibitem{BluCol}{\sc G.W. Bluman, J.D. Cole}, {\em Similarity Methods for Differential Equations}, Springer, Berlin, 1974.
%
\bibitem{Jost}{\sc J. Jost}, {\em Partial differential equations}, Springer, GTM 214, 2002.
%
\bibitem{Nad} {\sc M. Nadjafikhah}, {\em Lie Symmetries of Inviscid Burgers'
Equation}, Adv. appl. Clifford alg., DOI 10.1007/s00006-003-0000.
%
\bibitem{Olv1}{\sc P.J. Olver}, {\em Applications of Lie Groups to
Differential Equations}, Springer, New York, 1986.
%
\bibitem{Ovs}{\sc L.V. Ovsiannikov}, {\em Group Analysis of Differential
Equations}, Academic Press, New York, 1982.
%
\bibitem{PolZai} {\sc A.D. Polyanin}, and {\sc V.F. Zaitsev}, {\em Handbook of Nonlinear Partial Differential Equations}, Chapman \& Hall/CRC, Boca
Raton, 2004.
%
\bibitem{SechNa} {\sc R. Seshadri} and {\sc T.Y. Na}, {\em Group Invariance in Engineering Boundary Value
Problems}, Springer-Verlag, New York, 1985.
%
\bibitem{SekSha}{\sc T.R. Sekhar,V.D. Sharma}, {\em Similarity solutions for three dimensional Euler equations
using Lie group analysis}, Appl. Math. and Comp. 196 (2008)
147--157.
%
\end{thebibliography}
\end{document}